\documentclass[11pt]{article}
\newtheorem{thm}{Theorem}

\newtheorem{lemma}[thm]{Lemma}
\newtheorem{prop}[thm]{Proposition}

\newtheorem{ex}[thm]{Example}
\newtheorem{conj}[thm]{Conjecture}

\newcommand{\beq}[1]{\begin{equation}\label{#1}}
\newcommand{\enq}[0]{\end{equation}}

\newcommand{\qed}[0]{\begin{flushright} \rule{2mm}{3mm} \end{flushright}}

\newcommand{\C}[2]{{{#1}\choose{{#2}}}}
\newcommand{\ga}[0]{\alpha }
\newcommand{\gb}[0]{\beta }
\newcommand{\gc}[0]{\gamma }
\newcommand{\gd}[0]{\delta }

\newcommand{\gs}[0]{\sigma}
\newcommand{\gS}[0]{\Sigma}
\newcommand{\gz}[0]{\zeta}
\newcommand{\eps}[0]{\varepsilon }

\newcommand{\0}[0]{\emptyset}

\newcommand{\ra}[0]{\rightarrow}

\newcommand{\RR}[0]{\mbox{${\bf R}^+$}}
\newcommand{\Rr}[0]{\mbox{${\bf R}$}}

\newcommand{\E}[0]{{\sf E}}

\newcommand{\A}[0]{{\cal A}}
\newcommand{\B}[0]{{\cal B}}

\newcommand{\I}[0]{{\cal I}}
\newcommand{\m}[0]{{\cal M}}

\newcommand{\bn}[0]{\bigskip\noindent}
\newcommand{\mn}[0]{\medskip\noindent}

\newcommand{\sub}[0]{\subseteq}
\newcommand{\sm}[0]{\setminus}
\renewcommand{\dots}[0]{,\ldots,}

\newcommand{\uone}[0]{\underline{1}}

\newcommand{\cov}[0]{~\cdot\hspace{-0.09in}>}
\newcommand{\scov}[0]{~\cdot\hspace{-0.09in}\succ}

\newcommand{\da}[0]{\downarrow}
\newcommand{\ua}[0]{\uparrow}

\begin{document}

\renewcommand{\thefootnote}{\fnsymbol{footnote}}
\title{Negative correlation and log-concavity\footnotemark }

\author{
J. Kahn and M. Neiman\\
{Rutgers University } \\
{\footnotesize email: jkahn@math.rutgers.edu; neiman@math.rutgers.edu}
}
\date{}

\footnotetext{ * Supported by NSF grant DMS0701175.}

\maketitle


\begin{abstract}
We give counterexamples and a few positive results related to
several conjectures
of R. Pemantle \cite{Pem} and D. Wagner \cite{Wagner}
concerning negative
correlation and log-concavity properties for probability
measures and relations between them.
Most of the negative results have also been obtained, independently
but somewhat earlier, by Borcea {\em et al.} \cite{BBL}.
We also give short proofs of a pair of results from
\cite{Pem} and \cite{BBL}; prove that
``almost exchangeable" measures satisfy the ``Feder-Mihail" property,
thus providing
a ``non-obvious" example of a class of measures
for which this important property can be shown to hold;
and mention some further questions.
\end{abstract}

\section{Introduction}\label{Intro}

This paper is concerned with negative
correlation and log-concavity properties and relations
between them,
with much of our motivation provided by
\cite{Pem}
and \cite{Wagner}.
In particular, we give
counterexamples
to several conjectures from these papers,
and a positive answer to one question from \cite{choewagner03}.
While writing the present paper, we learned that
some of the conjectures disproved here were also recently
disproved in \cite{BBL}.
The present examples seem a little simpler and more natural,
and also show a little more, so are thought to still be
of interest.
We also give short proofs of some of the
results in \cite{BBL} and \cite{Pem},
and show that some versions of what was for us the
most interesting of the conjectures of \cite{Pem} fail even
for the natural example of competing urn measures,
the main point here being verification of a fairly strong negative
correlation property, ``conditional negative association,"
for such measures (proof of which will appear separately).
In this long introduction we summarize these developments
and some of the relevant open problems.

\bigskip
Given a finite set $S$,
denote by $\m=\m_S$ the set of probability measures
on $\Omega=\Omega_S=\{0,1\}^S$.
As a default we take
$S=[n] =\{1\dots n\}$ (which for us is simply
a generic $n$-set), using
$\Omega_n$ in place of $\Omega_{[n]}$.
%
We will occasionally identify
$\Omega$ with the Boolean algebra $2^{[n]}$
(the collection of subsets of $[n]$ ordered by inclusion)
in the natural way
(namely, identifying a set with its indicator).
Recall that an event $\A\sub \Omega$
is {\em increasing} (really, nondecreasing)
if $x\geq y\in\A$ implies $x\in\A$ (where we give $\Omega$ the product order),
and similarly for {\em decreasing}.
While our concern here is with negative dependence properties,
for perspective we first recall one or two points regarding
their better understood
positive counterparts.

\mn
{\em Positive correlation and association.}

Recall that events $\A,\B$ in a probability space
are \emph{positively
correlated}---we write $\A\uparrow \B$---if $\Pr(\A\B)\geq\Pr(\A)\Pr(\B)$.
The joint distribution of random variables $X_{1},...,X_{n}$---here
always $\{0,1\}$-valued---is said
to be \emph{positively associated} (PA)
if any two events both increasing
in the $X_{i}$'s are positively correlated. (This is easily seen
to be equivalent to the property that for any two increasing functions
$f,g$ of the $X_{i}$'s one has $\E(fg)\geq\E(f)\E(g)$.)

The seminal result here is {\em Harris' Inequality} \cite{Harris}, which says that
product measures are PA.
(The special case of uniform measure on $\Omega$
was rediscovered in \cite{Kleitman}, and in combinatorial circles has
often been called {\em Kleitman's Lemma}.)
The best known---and most useful---extension
of Harris' Inequality is the
{\em FKG Inequality} of Fortuin, Kasteleyn, and Ginibre
\cite{FKG}:
%
\begin{thm}\label{FKGthm}
If $\mu\in \m$ satisfies
\beq{PLC}
\mu(\eta)\mu(\tau)\leq\mu(\eta\wedge\tau)\mu(\eta\vee\tau)
~~~~~~~\forall\eta,\tau\in\Omega
\enq
(where $\wedge,\vee$ denote meet and join in the product order on
$\Omega$), then $\mu$ is PA.
\end{thm}
(Stronger still, and also very useful,
is the {\em Ahlswede-Daykin} or {\em ``Four Functions" Theorem}
\cite{AD}, whose statement we omit.)

The {\em positive lattice condition} (\ref{PLC}) (a.k.a.
the {\em FKG lattice condition} or
{\em log supermodularity}) is
equivalent to {\em conditional positive association},
the property that every
measure obtained from $\mu$ by conditioning on the values of some of the
variables is PA;
this follows easily from Theorem~\ref{FKGthm} and is a good
way to make sense of (\ref{PLC}).
One also says that $\mu$ with (\ref{PLC}) is an {\em FKG measure}.
See, e.g., \cite{Boll}, \cite{Grimmett}, \cite{LiggettIPS},
\cite{GHM}, \cite{DR}, \cite{FM} for a small sample of
applications of these notions in combinatorics,
probability, statistical mechanics, statistics and computer science.

\mn
{\em Negative association and related properties.}

While negative correlation
has the obvious meaning
($\mu(\A\B)\leq\mu(\A)\mu(\B)$, denoted $\A\downarrow \B$),
negative {\em association} requires a little care (for instance,
$\A\uparrow\A$ holds strictly for any $\A$ with
$\mu(\A)\not\in \{0,1\}$).
Say $i\in [n]$ {\em affects} event $\A$ if there are $\eta\in \A$
and $\tau\in\Omega\sm\A$ with $\eta_j=\tau_j$ for all $j\neq i$,
and write $\A\perp \B$ if no coordinate affects both $\A$ and $\B$.
Then $\mu\in \m$ is {\em negatively associated}
(or has {\em negative association}; we use ``NA" for either)
if
$\A\downarrow \B$
whenever $\A,\B$ are increasing and $\A\perp \B$.
We say that $\mu$ has {\em negative correlations} (or {\em is NC}) if
$\eta_i\downarrow \eta_j$
(that is, $\{\eta_i=1\}\downarrow \{\eta_j=1\}$)
whenever $i\neq j$.

Negative association turns out to be a much subtler
property
than PA.
Pemantle \cite{Pem} proposes
a number of questions regarding conditions related to NA, and
possible implications among them;
we sketch what we need from this,
and refer to \cite{Pem} for a more thorough discussion
(and more motivation).
The properties of interest for us are those obtained from NC
and NA by requiring closure under either conditioning or
imposition of external fields.
We first define these operations.

Unless specified otherwise, in this paper {\em conditioning} always means fixing
the values of some variables
(and this specification is always assumed to have
positive probability); thus a measure obtained from
$\mu\in \m$ by conditioning is one of the form
$\mu(\cdot|\eta_i=\xi_i ~\forall i\in I)$
for some $I\sub [n]$ and $\xi \in \{0,1\}^I$,
which we regard as a measure on $\Omega_{[n] \sm I}$.
(If we think of $\Omega$ as $2^{[n]}$, then
conditioning amounts to restricting our measure to
some interval $[J,K]$ of $2^{[n]}$ (and normalizing).)

For $W=(W_1\dots W_n)\in\Rr_+^n$ and $\mu\in\m$, define $W\circ \mu\in\m$ by
\beq{ExtFdef}
W\circ \mu(\eta) ~\propto~ \mu(\eta)\prod W_i^{\eta_i}
\enq
(meaning, as usual, that the left side is the right side
multiplied by the appropriate normalizing constant).
Borrowing Ising terminology, one says that $W\circ\mu$
is obtained from $\mu$ {\em by imposing the external field W}
(though to make this specialize correctly to the Ising model,
we should really take the ``field" to be $h$ given by
$h_i=\ln W_i$).
It will be convenient to allow $W_i=\infty$, which we interpret
as conditioning on $\{\eta_i=1\}$;
similarly we interpret $W_i=0$ as conditioning on $\{\eta_i=0\}$.


A third standard operation is projection:
the {\em projection} of $\mu$ on $J\sub [n]$ is the measure
$\mu'$ on $\{0,1\}^J$ obtained by integrating out
the variables of $[n]\sm J$; that is,
$$\mu'(\xi) =
\sum\{\mu(\eta):\eta\in \Omega, \eta_i=\xi_i ~\forall i\in J\}
~~~~(\xi\in \{0,1\}^J).$$

A basic motivation for much of \cite{Pem} was the
desire for a natural and robust notion (or notions)
of negative dependence, one measure of naturalness
(and also of usefulness) being invariance under some or
all of the preceding operations (and a few others that
we will not discuss here).  This leads in particular to
the following classes, which were alluded to above.

We say that $\mu\in \m$ is
\emph{conditionally negatively correlated} (CNC)
if every measure obtained from $\mu$ by conditioning is NC,
and NC+
if every measure obtained from $\mu$ by imposition of
an external field is NC.
\emph{Conditional negative association} (CNA)
and NA+ are defined analogously.
Of course NC+ and NA+ imply CNC and CNA respectively.
(This would be true even if we did not allow $W_i \in \{0,\infty\}$
in (\ref{ExtFdef}), since a limit of NC measures is
again NC, and similarly for NA; but there are properties
of interest---in particular the Feder-Mihail property below---for which
things go a little more smoothly with
the present convention.)

Note that Pemantle uses CNA+ where we use NA+,
but it is easy to see that the two notions coincide.
In general he uses ``+" for closure under
both projections and external fields, but for the properties
we are considering, this collapses to the definitions above:
it is easy to see
that all of the properties NC, CNC, NC+, NA, CNA, NA+ are
preserved by projections.

Following \cite{choewagner03}, \cite{Wagner},
we will also sometimes use the term \emph{Rayleigh} for NC+.
(The reference is to Rayleigh's monotonicity law for electric
networks;  see the second paragraph following Conjecture \ref{pef}
below or e.g. \cite{DS} or \cite{choewagner03}.)
We should also say a
little more about the relation between
our usage and that of \cite{Pem},
for which we need the
\emph{negative lattice condition} (NLC) for $\mu\in \m$:
\beq{NLC}
\mu(\eta)\mu(\tau)\geq\mu(\eta\wedge\tau)\mu(\eta\vee\tau)
~~~~~~~\forall\eta,\tau\in\Omega.
\enq
This is of course the analogue of (\ref{PLC}), but turns
out to be not nearly as useful, a crucial difference being
that, unlike (\ref{PLC}), it is not preserved by projections.
Following \cite{Pem}, we say that $\mu$ has the
\emph{hereditary negative lattice condition} (h-NLC) if
all projections of $\mu$ satisfy the NLC, and that
$\mu$ is \emph{h-NLC+} if every measure obtained from $\mu$
by imposition of an external field is h-NLC.
It is not hard to see that there are more names here than
properties:


\begin{prop} \label{prop:cnc_hNLC}

{\rm (a)} The properties CNC and h-NLC are equivalent.

\mn
{\rm (b)} The properties NC+ and h-NLC+ are equivalent.
\end{prop}
This has also been observed in \cite{BBL}
(see their Proposition 2.2 for (b) and Remark 2.2 for
a statement equivalent to (a)), so we will not prove it
here, but briefly:
(a) clearly implies (b);
h-NLC trivially implies CNC;
and the reverse implication follows easily
from the observation that the support of a
CNC measure is convex (i.e. $\mu(\eta),\mu(\tau)>0$
implies $\mu(\gs)>0$ whenever $\eta \leq \gs \leq \tau$),
proof of which is identical to
(e.g.) that of \cite[Theorem 4.2]{Wagner}.)

The next conjecture would be tremendously interesting.
\begin{conj}[\cite{Pem}]\label{Pem2and3}
{\rm (a)}  The properties CNC and CNA are equivalent.

\mn
{\rm (b)}  The properties NC+ and NA+ are equivalent.
\end{conj}
See \cite[Conjectures 2 and 3]{Pem}.
Note that in each case it is enough to show that the
first named property implies NA.
As shown in \cite{Pem}, CNA does not imply NC+;
we will later (Theorem \ref{Uthm})
see a ``naturally occurring"
example of this.
See also Conjecture~\ref{fmconj} below for one approach to proving
Conjecture \ref{Pem2and3}(b).

Recall that $\mu\in \m$ is {\em exchangeable} if it is invariant
under permutations of the coordinates (that is,
$\mu(\eta_{\gs(1)}\dots\eta_{\gs(n)}) = \mu(\eta_1\dots \eta_n)$
for any $\eta\in \Omega$ and permutation $\gs$ of $[n]$), or,
equivalently, if $\mu(\eta)$ depends only on $|\eta|:= \sum\eta_i$.
We say $\mu$ is {\em almost exchangeable} if it is
invariant under permutations of some subset of
$n-1$ of the variables.

Pemantle shows \cite[Theorem 2.7]{Pem} that for symmetric measures
the properties CNC, NC+, CNA and NA+ are equivalent,
while \cite{BBL} proves Conjecture \ref{Pem2and3}
for almost exchangeable measures:

\begin{thm}[\cite{BBL}, Corollary 6.6] \label{AE_equiv}
For almost exchangeable measures
\\{\rm (a)} the properties CNC and CNA are equivalent, and
\\{\rm (b)} the properties NC+ and NA+ are equivalent.
\end{thm}
In Section \ref{section:ae} we give quick proofs of both these results.
(Note that, in contrast to the exchangeable case,
CNA and NC+ are not equivalent for almost exchangeable measures;
see Theorem \ref{Uthm} and Example \ref{Uex2}.)
It may be worth noting that, despite its apparent simplicity,
the class of almost exchangeable measures is considerably richer
than the class of exchangeable measures;
in particular, the examples proving Theorem \ref{Tctrex} below
(and also those of \cite{BBL}) are almost exchangeable.

\mn
{\em Log-concavity}

Recall that a sequence $a=(a_0\dots a_n)$ of real numbers
is {\em unimodal} if there is some $k\in \{0\dots n\}$ for which
$a_0 \leq a_1 \leq \cdots \leq a_k \geq \cdots \geq a_n $,
and is {\em log-concave} (LC) if
$a_i^2\geq a_{i-1}a_{i+1}$ for $1\leq i\leq n-1$.
Of course a nonnegative LC sequence with no internal zeros is unimodal.
Following \cite{Pem} we say that $a$ (as above)
is {\em ultra-log-concave} (ULC) if
the sequence $(a_i/\C{n}{i})_{i=0}^n$ is log-concave
and has no internal zeros.

We also say that
$\mu\in\m$ is ULC if
its {\em rank sequence}, $(\mu(|\eta|=i))_{i=0}^n$, is ULC.
We define ``$\mu$ is LC" and ``$\mu$ is unimodal" similarly,
except that for the former we add the stipulation that
the rank sequence has no internal zeros.
(It would be convenient to also make this a requirement for
``LC" for sequences, but we politely adhere to the standard
definition.)

Pemantle shows \cite[Theorem 2.7]{Pem} that for exchangeable
measures, ULC coincides with CNC, NC+, CNA and NA+.
He conjectures (see his Conjecture 4)
that each of the latter properties
implies ULC for general $\mu$; more precisely, this is a set
of four conjectures, the weakest of which is
the one with the strongest hypothesis:

\begin{conj} \label{Pem4}
NA+ implies ULC.
\end{conj}
(He also conjectures that NA implies ULC, but, as noted in
\cite{BBL}, this is easily seen to be incorrect, even for
exchangeable measures.)

One of the stronger versions of Pemantle's conjecture---that
the Rayleigh property
(i.e. NC+) implies ULC---was separately
proposed by Wagner in \cite{Wagner}, where it was called
the ``Big Conjecture."
(The overlap seems due to the failure in \cite{Pem}, \cite{Wagner}
to notice Proposition \ref{prop:cnc_hNLC}(b).)
We will say more about Wagner's motivation below.
In Section \ref{Ctrx} we
give a family of examples that disproves all these conjectures and more:
\begin{thm}\label{Tctrex}
Conjecture \ref{Pem4} is false;
in fact NA+ does not even imply unimodality.
\end{thm}

As mentioned earlier, we recently learned that
the first part of this was discovered a little
earlier by Borcea {\em et al.} \cite{BBL}.
The present examples are slightly smaller
(12 variables as opposed to 20 for violation of ULC) and simpler,
and also disprove more, as the
example of \cite{BBL} is LC.

The examples for Theorem \ref{Tctrex} also turn out to
disprove Conjectures 8 and 9 of \cite{Pem};
again, the first of these is also disproved by the
example of \cite{BBL}.
Statements of these conjectures are deferred to Section \ref{Ctrx}.

\medskip
A more particular notion than ULC, from \cite{choewagner03},
is as follows.
For positive integer $m$, $\mu\in \m$ is
said to have the property LC[$m$] if, for every $S \sub [n]$
of size at most $m$, every measure obtained from $\mu$ by
imposing an external field and then projecting on $S$ is ULC.
Choe and Wagner \cite[Theorem 4.6]{choewagner03} show that
the three properties NC+, LC[2], and LC[3] are equivalent.
(Strictly speaking,
\cite{choewagner03} is confined to a smaller class
of $\mu$'s, but the proof is valid in the present generality.)
They ask whether NC+ implies LC[4].
(Of course,
since projections preserve NC+, Wagner's conjecture above
would say that NC+ implies LC[$m$] for every $m$.)
The next result is proved in Section \ref{sec:rlc5}.
\begin{thm} \label{thm:rlc5}
NC+ implies LC[5].
\end{thm}
Thus NC+ also
implies LC[4], whereas the examples for Theorem \ref{Tctrex}
will show
that
\beq{lc12}
\mbox{NC+ does not imply LC[12].}
\enq
We don't know what happens between 5 and 12.
Of course Theorem \ref{thm:rlc5}
is now less interesting than formerly, when it was thought
to be a step in the direction of Conjecture \ref{Pem4}.

\medskip
Let us pause here to mention a strengthening of ULC.
For $\mu\in \m$ set
$$\ga_i(\mu) =
\C{n}{i}^{-1}
\sum_{\eta\in \Omega, |\eta|=i} \mu(\eta)\mu(\uone-\eta)$$
(where $\uone=(1 \dots 1)$).
Say that $\mu\in \m_{2k}$ has the {\em antipodal pairs property}
(APP) if
$\ga_k(\mu)\geq \ga_{k-1}(\mu),$
and that $\mu\in \m_n$ has the
{\em conditional antipodal pairs property} (CAPP)
if any measure obtained from $\mu$ by conditioning on
the values of some $n-2k$ variables (for some $k$) has the APP.
(Note these properties are not affected by external fields.)

\begin{thm}\label{TAPP}
The CAPP implies ULC.
\end{thm}
This is reminiscent of an observation of T. Dowling \cite{Dowling};
see the paragraph following Conjecture \ref{CAPP}.
As pointed out to us by David Wagner \cite{WagnerPC},
Theorem 4.3 of \cite{Wagner2}, which was motivated by
\cite{Dowling}, feels similar to Theorem \ref{TAPP}.
In fact, Theorem \ref{TAPP} implies strengthenings of Theorem 4.3 and a few other results in \cite{Wagner2}; details of this connection will appear elsewhere \cite{KN3}.

Though we will not do so here, we can strengthen
Theorem \ref{thm:rlc5} to its antipodal pairs
version (gotten by replacing
ULC by CAPP in the definition of LC[$m$]).
The proof of Theorem \ref{TAPP}
requires some work, depending, {\em inter alia,} on
Delsarte's inequalities \cite{Delsarte}, and
will appear elsewhere \cite{KN3}.

\mn
{\em Mason's Conjecture}

Here we want to say a little about the motivation for
Wagner's (``big") conjecture and mention a few related questions.
For this discussion we regard a {\em matroid} as a collection
$\I$ of independent sets, subsets of some ground set $E$.
We will not go into matroid definitions; see e.g.
\cite{Welshbook} or \cite{Oxley}.
Prototypes are the collection of (edge sets of)
forests of a graph (with edge set $E$)---this is a
{\em graphic} matroid---and (as it turns out, more generally)
the collection of linearly independent subsets of
some finite subset $E$ of some (not necessarily finite)
vector space.
For present purposes not too much is lost
by thinking only of graphic matroids.

We are interested in the
{\em independence numbers} of a matroid $\I$, that is, the numbers
$$a_k = a_k(\I) = |\{I\in\I:|I|=k\}| ~~~ k=0\dots n,$$
concerning which we have a celebrated conjecture of
J. Mason:

\begin{conj}[\cite{Mason}]\label{Cmason}
For any matroid $\I$ on a ground set of size n,
the sequence $a=a(\I)=(a_0\dots a_n)$ is ULC.
\end{conj}
(Note that $a$ will typically end with some zeros,
and also that in the graphic case $n$ counts edges, not vertices.)
Of course one can relax Conjecture \ref{Cmason}
by asking for LC or unimodality
in place of ULC.  In fact unimodality,
first suggested by Welsh \cite{Welsh},
was the original conjecture
in this direction, and
even this, even for graphic matroids, remains open.
(See \cite{Stanley}
or \cite{Brenti}
for much more on
log-concavity in combinatorial settings.)

From the present viewpoint, Mason's Conjecture asks for
ultra-log-concavity of uniform measure on $\I$ (regarded
in the usual way as a subset of $\{0,1\}^E$).
In case $\I$ is graphic such a measure is a
{\em uniform spanning forest} (USF) measure
(``spanning" because we think of a member of $\I$ as a
subgraph that includes all vertices).
These measures are also very interesting from a
correlation standpoint.

\begin{conj}\label{pef}
USF measures are Rayleigh.
\end{conj}
This natural guess was perhaps first proposed in \cite{AN}
(which was circulated in the combinatorial community
as early as 1993, but took a while to get to press).
It is also, for example, Conjecture 5.11.2 in \cite{Wagner}.
(The statement in \cite{AN}
is (in present language)
that USF measures are NC, but it is not hard to see that this is
equivalent.)

As essentially shown by Brooks et al. \cite{BSST},
the analogue of Conjecture \ref{pef}
for uniform measure on
the spanning {\em trees} of a (finite) graph
amounts to Rayleigh's monotonicity law for electric networks
(again, see \cite{DS}).  This was extended by
Feder and Mihail (\cite{FM}, to which we will return shortly)
to say that such measures
are in fact NA+ (more precisely, this is what their proof gives).

For use below, let us call a measure obtained from a
USF measure by imposition of an external field---equivalently,
a measure $\mu$ on the spanning forests of
some finite graph $G$ with, for some $W:E(G)\ra \RR$,
$\mu(F)\propto \prod_{e\in F} W(e)$---a {\em weighted spanning
forest (weighted SF)}
measure, and define {\em weighted spanning tree (WST)} measures
and {\em weighted matroid measures}
(replace ``forest" by
``independent set") similarly.
(We avoid ``WSF" since it means {\em wired} SF; see e.g. \cite{LP}.)

One should note that,
while the intuition for Conjecture~\ref{pef} may seem
clear---presence of a given edge $e$ makes it easier for
a second edge $f$ to complete a
cycle---this may be misleading, since
the same intuition applies to
uniform measure on
the independent sets of a general matroid,
for which NC need {\em not} hold
(as can be derived from an example of
Seymour and Welsh \cite{SW}).
Some evidence for Conjecture \ref{pef},
and its analogue for spanning connected subgraphs,
is given in
\cite{GW}.
Also worth mentioning here---though without definitions; see
\cite{RCM}---is the following far-reaching extension of
Conjecture \ref{pef},
which has been ``in the air" for a while
(e.g. \cite{Pem}, \cite{Gsur}).
\begin{conj}\label{CRCM}
Any random cluster measure with $q<1$ is NA+.
\end{conj}
(Equivalently, such measures are NA.)
Limiting cases include the aforementioned uniform measures on
forests, spanning trees and connected subgraphs of a graph;
again see \cite{RCM}.
Conjecture \ref{CRCM} with NC+ in place of NA+ is proved for
series-parallel graphs (part of a more general matroid statement) in
\cite{Wagner} (see Example 5.1 and Theorem 5.8(d)).

\medskip
Of course Wagner's ``Big Conjecture," if true, would
have implied
Mason's Conjecture for any class of matroids for which
one could establish the Rayleigh property (meaning, of course,
for uniform measure on independent sets).
Conjecture \ref{pef} says that graphic matroids
should be such a class, and
Wagner \cite[Conj. 5.11]{Wagner}
suggests a sequence of strengthenings of this.
(Mason's Conjecture also
partly motivated Conjecture \ref{pef} in \cite{AN},
though at the time
the connection was not much more than a feeling that the
issues underlying
the two were similar.)

\medskip
Let us also mention that, as far as we know,
the following strengthening of Mason's Conjecture
could
be true.

\begin{conj}\label{CAPP}
Weighted matroid measures have the CAPP.
\end{conj}
(Of course it is enough to prove APP.)
We can now fill in our earlier allusion to \cite{Dowling}:
Dowling observed that the LC version of
Mason's Conjecture would follow from
the assertion that, for any matroid $\I$ on a ground set $E$
of size $2k$, the number of ordered partitions $E=I\cup J$
with $I,J\in \I$ and $|I|=|J|=k$ is at least as large
as the number with $|I|=k-1$ (and $|J|=k+1$),
and showed this is true for $k\leq 7$.
In fact, Dowling's proof also shows that, for $k \leq 5$, uniform measure on the independent sets of any matroid on a ground set of size $2k$ has the APP.
This gives Conjecture \ref{CAPP} for matroids on ground sets of
size at most $11$ and, via
(a slight generalization of) Theorem \ref{TAPP}, proves that,
for any matroid on a set of size $n$, the sequence $\{a_i\}$
of independence numbers is ``ULC up to 6,'' meaning
$(a_k/\C{n}{k})^2 \geq (a_{k-1}/\C{n}{k-1}) (a_{k+1}/\C{n}{k+1})$
for $k\leq 5$;
see \cite{KN3} for details.
This is a small improvement on the best that seems to have been known
previously, namely that $\{a_i\}$ is ULC up to 4,
which was shown by
Hamidoune and Sala\"un \cite{HamidouneSalaun}.

\mn
{\em Feder-Mihail}

Say $\mu\in\m$ has
the {\em Feder-Mihail property}
(or {\em is FM}) if

\mn
$~~~~~~~~~~~~~~~~~~~~~~~~$
{\em for any increasing $\A\sub \Omega$,
$~\{\eta_i=1\}\uparrow \A$ for
some $i\in [n]$,}

\mn
and extend this to CFM and FM+ in the usual way.
(Of course FM+ trivially implies CFM, but note that this
implication requires that we explicitly include conditioning
in our definition of ``+''
(i.e. we allow $W_i \in \{0,\infty\}$ in (\ref{ExtFdef})),
since e.g. there are situations where $W \circ \mu$
is FM for all $W$ with positive entries
but $\mu(\cdot|\eta_1=1)$ is not FM.)
The following simple but powerful observation is
essentially from \cite{FM}, though given there only in
a special case.

\begin{thm}\label{TFM}
{\rm (a)}  If $\mu\in \m$ is both CNC and CFM then it is CNA.

\mn
{\rm (b)}  If $\mu\in \m$ is both NC+ and FM+ then it is NA+.

\end{thm}
(A statement equivalent to (a) is proved in
\cite[Theorem 1.3]{Pem}, and (b) follows easily from (a).)

Given the power of Theorem \ref{TFM}, it would be useful
to identify situations where the FM property
holds.
This is trivially the case for $\mu$ concentrated on a level
(that is, $\{\eta\in\Omega: |\eta|=k\}$ for some $k$; see e.g. \cite[Corollary 3.2]{DJR}),
e.g. for WST measures (a key to \cite{FM}),
and it is fairly easy to show that exchangeable measures and,
more generally, ``rescalings" of product measures,
satisfy the stronger ``normalized matching property"
(NMP; see Section \ref{section:ae} for definitions).
But in general FM seems hard to establish,
and indeed we are not aware of any interesting classes of
non-NMP measures that are known to be FM.
Thus the following result, which is proved in
Section \ref{section:ae}, may be of some interest.

\begin{thm}\label{TAE}
Almost exchangeable measures are FM+.
\end{thm}
Note that this combined with Theorem \ref{TFM} gives
Theorem \ref{AE_equiv}.
(This is not quite the ``quick" proof
of Theorem \ref{AE_equiv} promised earlier, since
Theorem \ref{TAE} requires some effort; but, as observed in
Section \ref{section:ae}, FM (resp. FM+)
for almost exchangeable measures that are also NC (resp. NC+) is much easier.)

Despite the (apparent) difficulty of proving FM, the property
seems to tend to hold for measures not deliberately
constructed to violate it.
In particular
we would like to propose, perhaps a bit optimistically,
the following possibilities.
%
\begin{conj} \label{fmconj}
The Feder-Mihail property holds for
\\
{\rm (a)}  Rayleigh measures,
\\
{\rm (b)} weighted SF measures, and (more generally)
\\
{\rm (c)} weighted matroid measures.
\end{conj}
Note that, in view of Theorem \ref{TFM}, (a) would imply
Conjecture \ref{Pem2and3}(b)
(the corresponding approach to Conjecture \ref{Pem2and3}(a) fails
because CNC measures need not be FM),
while (b)
together with Conjecture~\ref{pef} would
say that USF measures are NA+.
(Extending this to matroids via (c) fails because
Conjecture~\ref{pef} does.)

\mn
{\em Competing urns}

One of the principal motivating examples for \cite{Pem}
is ``competing urns."
Here we have $m$ balls which are thrown independently
into urns $1\dots n$ according to some (common) distribution,
let $X_i$ be the indicator for occupation of urn $i$,
and consider the corresponding measure $\mu$ on $\Omega$
(that is, the law of $(X_1\dots X_n)$).
(Formally we may take a random $\gs:[m]\ra [n]$
with the $\gs(i)$'s i.i.d., and
$X_i={\bf 1}_{\{\gs^{-1}(i)\neq \0\}}$.)
Call a $\mu$ of this type a {\em competing urn measure}.

The competing urns model is
explored in some detail by Dubhashi and Ranjan \cite{DR}, who
in particular prove that competing urn measures are NA.
Another proof of this is given in \cite{Pem}.

More generally one may consider thresholds $a_1\dots a_n$,
with $X_i $ the indicator of $\{|\gs^{-1}(i)|\geq a_i\}$---for lack
of a better name, we then call the law of $(X_1\dots X_n)$ an
{\em extended} competing urn measure---and
the arguments of \cite{DR} and \cite{Pem} apply to show that such
measures are again NA.
(Actually \cite{DR} proves the stronger statement that
the (law of the) random variables
$X_i(j)={\bf 1}_{\{\gs(i)=j\}}$ is NA.)

Here again, we have a suggestion from \cite{Pem}:

\begin{conj}\label{Pem5}
Any competing urn measure is ULC.
\end{conj}
In fact
Pemantle
conjectures something more general that
we will not state,
since, unfortunately,
Conjecture \ref{Pem5} is
not true.

\begin{prop}\label{Uctrex}
Competing urn measures need not be LC.
\end{prop}
(See Example \ref{Uex1}.)
It is also not hard to give examples of non-Rayleigh
competing urn measures
(Example \ref{Uex2}).
On the other hand,
one may argue that in the context of competing urns,
external fields are less natural
than
conditioning;
at any rate, it does turn out that even
extended competing urn measures are nicely behaved under conditioning:
\begin{thm}\label{Uthm}
Extended competing urn measures are CNA.
\end{thm}
Thus, as mentioned earlier, we have a natural class of
measures for which CNA does not imply Rayleigh,
and, combining with Proposition \ref{Uctrex},
a (natural) counterexample to the strengthening of
Conjecture \ref{Pem4} obtained by replacing NA+ by CNA
(this is
again one of the versions of Conjecture 4 of \cite{Pem}).
%

While
Theorem \ref{Uthm} seems like it ought to be easy,
we do not know a simple argument, even in the ordinary
(non-extended) case,
and, to keep this paper from getting too long,
will give the proof elsewhere \cite{KN2}.
That it might not be completely straightforward is
suggested by our inability to decide whether it is still true
(again, even in the ordinary case) if we
drop the requirement
that the balls be identically distributed.
The aforementioned NA for $X_i(j)$'s proved in \cite{DR}
is true in this generality, which gives NA for these more general
urn measures.
(The argument of \cite{Pem} does not work with nonidentical balls.)

\section{Exchangeable and almost exchangeable measures}\label{section:ae}

We need a few more definitions.
We extend the definitions of exchangeable and almost exchangeable measures
(given following Conjecture \ref{Pem2and3}) to
general functions on $\Omega$ in the obvious way
($f:\Omega \to \Rr$ is {\em almost exchangeable} if it is invariant under
permutations of some subset of $n-1$ of the variables and {\em exchangeable}
if it is invariant under permutations of all the variables).
We also extend our notation for positive and negative correlation to functions:
for $f,g:\Omega \to \Rr$, we write $f \ua g$ if $\E(fg) \geq \E(f) \E(g)$
(and similarly for $f \da g$); we will also write, e.g.,
$\A \ua f$ for ${\bf 1}_{\A} \ua f$.
A stronger statement is that
$\A$ is {\em stochastically increasing in $f$}, that is, that
$\Pr(\A|f=t)$ is increasing in $t$,
where we restrict to values of $t$ for which $\Pr(f=t)$
is positive.
(Note that our use of the notation $\A \ua f$ differs from that in \cite{Pem}.)
Following \cite{DJR}, for a function $f:\Omega \to \Rr$ and a measure $\mu \in \m$, we
say that $i\in [n]$
is a {\em variable of positive influence} for the pair $(f,\mu)$
(or $(\A,\mu)$ if $f={\bf 1}_{\A}$) if $\eta_i \ua f$.
Thus the FM property for $\mu$ says that for every increasing $\A$ there is a
variable of positive influence for $(\A,\mu)$.
In \cite{DJR}, Dubhashi {\em et al.} prove the (easy) result that $(f,\mu)$ has a
variable of positive influence if $f$ is increasing and at least one of $f$, $\mu$ is
exchangeable, and ask for other classes of function-measure pairs having variables of
positive influence.
Here we prove the following result, which, as we will see shortly,
implies Theorem \ref{TAE}.
\begin{thm}\label{aepi}
If there is a variable $l$ for which
\beq{incr_conds}
\E[f|\eta_l=j,\sum_{i \neq l} \eta_i=k] \textrm{ is increasing in } j \textrm{ and } k
\enq
(for pairs $(j,k)$ for which the conditioning event has positive probability), then $(f,\mu)$ has a variable of positive influence.
In particular, if $f$ is increasing and almost exchangeable then $(f,\mu)$ has a variable of positive influence.
\end{thm}

As observed earlier, Theorem \ref{AE_equiv} is an immediate consequence of
Theorems \ref{TFM} and \ref{TAE}; but it does not require the
full strength of
Theorem \ref{TAE}, and before proving Theorem \ref{aepi} we will give an easier
argument, together with a quick proof of the following result from \cite{Pem}.
\begin{thm}\label{E_equiv}
For exchangeable measures the properties CNC, CNA, NC+, NA+, and ULC are equivalent.
\end{thm}

For these arguments and the derivation of Theorem \ref{TAE}
we need a little background.
We first recall Chebyshev's Inequality, here stated in our terminology.
\begin{prop}\label{ChebIneq}
Any probability measure on a totally ordered set is PA.
\end{prop}
(Where, of course, PA is as for measures on $\{0,1\}^S$:
$f \ua g$ for any two increasing functions $f,g$.)

Recall that for probability measures $\mu$ and $\nu$, $\mu$
\emph{stochastically dominates} $\nu$ (written $\mu \succeq \nu$) if
$\mu(\A) \geq \nu(\A)$ for every increasing $\A$.
Writing $\mu_k$ for the conditional measure $\mu(\cdot|\sum \eta_i=k)$
(defined only when $\mu(\sum \eta_i=k)>0$), we say that $\mu$ has the
{\em normalized matching property} (NMP) if $\mu_{l} \succeq \mu_k$ whenever
$l \geq k$ and both conditional measures are defined.
(This generalizes the usual definition, for which see e.g. \cite{Anderson}.)
The NMP is equivalent to the property that every increasing event $\A$ is
stochastically increasing in $\sum \eta_i$,
which
implies (easily and directly, by Proposition \ref{ChebIneq}, or
essentially by Proposition 1.2 in \cite{Pem})
that
$\A \ua \sum \eta_i$,
and thus
(since expectation is linear) that $\A \ua \eta_i$ for some $i$.
This discussion gives the next observation.

\begin{prop}\label{nmp_fm}
The NMP implies FM.
\end{prop}
Conjecture 8 of \cite{Pem} says that NA+ implies NMP, but this is false; see Conjecture \ref{conj:Pem8} and Theorem \ref{thm:ctrex_Pem89} in Section \ref{Ctrx}.
Given $\mu \in \m$ and nonnegative sequence $a=(a_i)_{i=0}^n$, the {\em generalized rank rescaling of $\mu$ by $a$} is the measure $a \otimes \mu \in \m$ with
$$a \otimes \mu(\eta) \propto a_{|\eta|} \mu(\eta).$$
(To be precise, we only make this definition when the right side
is not identically zero.)
This generalizes the rank rescaling operation in \cite{Pem}, which required that $a$ be LC with no internal zeros.
Observe that (since $\mu_k=(a \otimes \mu)_k$ whenever $a_k > 0$) generalized rank rescalings preserve the NMP.

\begin{lemma} \label{rr_prod_nmp}
Product measures (and, consequently, generalized rank rescalings of product measures) have the NMP.
\end{lemma}
(A proof is sketched in \cite[Section 4.2]{DJR}. More generally, any product of LC measures with the NMP again has these properties;
this was proved by Harper \cite{Harper}, and can also
essentially be gotten from a combinatorial version proved in
\cite{Hsieh-Kleitman}, \cite{Anderson}.)
For the proof of Theorem \ref{AE_equiv} we also need
the following standard observation,
an easy consequence of Proposition \ref{ChebIneq}.
\begin{lemma} \label{little_lemma}
Let $f,g:\Omega \to \Rr$, and suppose for some event $\B$
\\{\rm (i)} each of $f$, $g$ is positively correlated with $\B$, and
\\{\rm (ii)} $f$ and $g$ are conditionally positively correlated given each of $\B$, $\Omega \sm \B$.
\\Then $f$ and $g$ are positively correlated.
\end{lemma}


\mn{\em Proof of Theorem \ref{TAE}.}
Let $\mu' \in \m$ be invariant under permutations of the
variables $1 \dots n-1$,
and
$\mu=W \circ \mu'$
for some $W \in \Rr_+^n$.
We may assume all $W_i$ are finite and strictly positive,
since otherwise we can reduce the number of variables
(note that any measure gotten from an almost exchangeable measure
by conditioning is again almost exchangeable).
Then, with $\nu \in \m_{n-1}$ the product measure satisfying
$$\nu(\eta) \propto \prod_{i \in [n-1]} W_i^{\eta_i},$$
we have
$$\mu(\cdot|\eta_n=0, \sum_{i \in [n-1]} \eta_i=k) = \mu(\cdot|\eta_n=1, \sum_{i \in [n-1]} \eta_i=k) = \nu_k,$$
so (by Lemma \ref{rr_prod_nmp})
$f={\bf 1}_{\A}$ satisfies
(\ref{incr_conds}) with $l=n$ for every increasing event $\A$. \qed

\mn{\em Proof of Theorem \ref{AE_equiv}.} By Theorem \ref{TFM},
it suffices to show that
{\em every NC measure that can be gotten by applying an
external field to an almost exchangeable measure is FM}.
Let $\mu$ be such a measure, obtained by imposing an external field on a
measure invariant under permutations of coordinates $\{2 \dots n\}$, and let
$\A$ be an increasing event.
We should show that $\A \ua \eta_i$ for some $i \in [n]$.
We may assume that all coordinates of the external field are finite
and strictly positive (or we can reduce the number of variables),
and that $\A \da \eta_1$ (or we are done).
For $j \in \{0,1\}$, the conditional measure $\mu(\cdot|\eta_1=j)$ is a generalized rank rescaling of a product measure, so by Lemma \ref{rr_prod_nmp} and the discussion in the paragraph preceding Proposition \ref{nmp_fm} we have $\A \ua f := \sum_{i \neq 1} \eta_i$ conditionally given either of the events $\{\eta_1=0\}$, $\{\eta_1=1\}$.
Since $\mu$ is NC, we have $\eta_1 \da \eta_i$ for all $i \in \{2 \dots n\}$,
so that $\eta_1 \da f$.
But then
applying Lemma \ref{little_lemma} with $g={\bf 1}_{\A}$ and $\B=\{\eta_1=0\}$
gives $\A \ua f$, whence
$\A \ua \eta_i$ for some $i \in \{2 \dots n\}$. \qed

\mn{\em Proof of Theorem \ref{E_equiv}.} It suffices to show that CNC implies ULC, and that ULC implies NA+.
The first implication is easy: CNC implies NLC (cf. Proposition \ref{prop:cnc_hNLC}), which for exchangeable measures is equivalent to ULC.
For the second implication, our main point is that we can eliminate much of the
work in \cite{Pem} by observing that exchangeable measures are FM+ (by Lemma
\ref{rr_prod_nmp} and Proposition \ref{nmp_fm}; of course this is also an instance of
Theorem \ref{TAE}, but not one that requires the less trivial Theorem \ref{aepi}),
so that by Theorem \ref{TFM}(b) it is enough to show
that ULC implies NC+.
This is a special case of the observation, proved in Lemma 2.8 of \cite{Pem}, that a
measure obtained from an exchangeable ULC measure by imposing an external field that
is identically $1$ on $J \sub [n]$, followed by projection on $J$, is exchangeable and
ULC. \qed

\mn{\em Proof of Theorem \ref{aepi}.}
First observe that if $f:\Omega_n \to \Rr$ is increasing
and invariant under permutations of coordinates in $[n] \sm \{l\}$
then (\ref{incr_conds}) is satisfied for every $\mu \in \m$, so the
last part of Theorem \ref{aepi} follows from the first.

To prove the first part, suppose, without loss of generality, that $\mu \in \m$ and
$f:\Omega_n \to \Rr$ satisfy (\ref{incr_conds}) with $l=1$, and set $h(\eta)=\sum_{i \neq 1} \eta_i$ ($\eta\in \Omega_n$).
It suffices to show that either
\beq{eitheror}
\mbox{$f \ua \eta_1 \ \ $ or $ \ \ f \ua h$.}
\enq

For $i \in \{0 \dots n-1\}$, let
$$\ga_i=\mu(h=i,\eta_1=1) \ \ \ \textrm{ and } \ \ \ \gb_i=\mu(h=i,\eta_1=0).$$
Choose increasing, nonnegative sequences
$\gc = (\gc_0\dots \gc_{n-1})$ and $\gd = (\gd_0\dots \gd_{n-1})$
such that
\beq{gcgd}
\gc_i=\E[f|h=i,\eta_1=1] \ \ \ \textrm{ and } \ \ \ \gd_i=\E[f|h=i,\eta_1=0]
\enq
whenever the conditioning events have positive probability
and $\gc_i \geq \gd_j$ whenever $i \geq j$.
(Existence of $\gc,\gd$ is guaranteed by (\ref{incr_conds}).
This extension to values not given by (\ref{gcgd}) is convenient,
but not really necessary, as these values play no role; see
(\ref{pc_fh}) and Lemma \ref{ETS_for_aepi}.)

Assume $f \da \eta_1$, i.e.,
$$\frac{\sum \ga_i \gc_i}{\sum \ga_i} \leq \frac{\sum \gb_i \gd_i}{\sum \gb_i}$$
(all sums in this proof are over $\{0 \dots n-1\}$ unless otherwise specified).
We want to show $\E(fh)\geq \E(f) \E(h)$,
that is,
\beq{pc_fh}
\sum (\ga_i \gc_i + \gb_i \gd_i)  \sum i (\ga_i + \gb_i) \leq
\sum i(\ga_i \gc_i + \gb_i \gd_i) \sum (\ga_i + \gb_i) .
\enq
(Of course the last sum is 1.)
This will follow from
$$
 \sum \ga_i \gc_i   \sum i \ga_i  \leq
 \sum i \ga_i \gc_i   \sum \ga_i,
$$
$$
 \sum \gb_i \gd_i   \sum i \gb_i  \leq  \sum i \gb_i \gd_i   \sum \gb_i,
$$
and
\beq{ineqiii}
 \sum i \ga_i   \sum \gb_i \gd_i  +  \sum i \gb_i   \sum \ga_i \gc_i
 \leq  \sum i \ga_i \gc_i   \sum \gb_i  +  \sum i \gb_i \gd_i  \sum \ga_i.
\enq
The first two of these
are instances of Proposition \ref{ChebIneq} (since $\gc$ and $\gd$ are increasing),
so it suffices to prove the following.
\begin{lemma}\label{ETS_for_aepi}
Let $\ga=(\ga_i)_{i=0}^{n-1}$, $\gb=(\gb_i)_{i=0}^{n-1}$, $\gc=(\gc_i)_{i=0}^{n-1}$, and $\gd=(\gd_i)_{i=0}^{n-1}$ be nonnegative sequences
(with neither of $\ga$, $\gb$ identically zero).
If $\gc$ and $\gd$ are increasing, $\gc_i \geq \gd_j$
whenever $i \geq j$, and $(\sum \ga_i \gc_i)/(\sum \ga_i) \leq
(\sum \gb_i \gd_i)/(\sum \gb_i)$, then (\ref{ineqiii}) holds.
\end{lemma}
{\em Proof.} Since scaling $\ga$, $\gb$ affects neither our hypotheses
nor (\ref{ineqiii}), we may assume $\sum \ga_i = \sum \gb_i$.
It suffices to show
\beq{splitineq}
 \sum_{i \geq s} \ga_i   \sum \gb_i \gd_i  +  \sum_{i \geq s} \gb_i  \sum \ga_i \gc_i \leq  \sum_{i \geq s} \ga_i \gc_i   \sum \gb_i  + \sum_{i \geq s} \gb_i \gd_i   \sum \ga_i
\enq
for $s \in [n-1]$ (since summing (\ref{splitineq}) over $s$
yields (\ref{ineqiii})).

Fix $s \in [n-1]$.
Obviously,
\beq{mono}
\mbox{if (\ref{splitineq}) is true, then it remains true when any $\gd_i$ with $i \geq s$ is increased.}
\enq
We define $\gd'=(\gd'_i)_{i=0}^{n-1}$ by
$$\gd'_i=\left\{ \begin{array}{ll}
\gd_i & \textrm{ if } i < s \\
\gd_s & \textrm{ if } i \geq s
\end{array} \right.$$
and consider two cases.

\mn
\underline{Case 1}:
$\sum \ga_i \gc_i> \sum \gb_i \gd'_i$.
Then there is an increasing sequence $\gd''=(\gd''_i)_{i=0}^{n-1}$ with
$\gd'_i \leq \gd''_i \leq \gd_i$ for all $i$ and
$\sum \ga_i \gc_i=\sum \gb_i \gd''_i$ (note $\sum \ga_i \gc_i \leq \sum \gb_i \gd_i$, since we normalized to $\sum \ga_i = \sum \gb_i$).
Since $\gc$ and $\gd''$ are increasing, we have
$$ \sum_{i \geq s} \ga_i   \sum \ga_i \gc_i
\leq  \sum_{i \geq s} \ga_i \gc_i   \sum \ga_i $$
and
$$ \sum_{i \geq s} \gb_i   \sum \gb_i \gd''_i  \leq  \sum_{i \geq s} \gb_i \gd''_i \sum \gb_i .$$
This yields (\ref{splitineq}) with $\gd$ replaced by $\gd''$, and  then  (\ref{splitineq}) (for $\gd$) follows from (\ref{mono}).

\mn
\underline{Case 2}:
$\sum \ga_i \gc_i \leq \sum \gb_i \gd'_i$.
By (\ref{mono}), it suffices to prove (\ref{splitineq}) with $\gd$ replaced by $\gd'$; this is a straightforward computation:
\begin{eqnarray*}
\sum_{i \geq s} \ga_i   \sum \gb_i \gd'_i  +  \sum_{i \geq s} \gb_i  \sum \ga_i \gc_i
& \leq &  \sum_{i \geq s} (\ga_i + \gb_i) \sum \gb_j \gd'_j  \\
& = & \sum_{i \geq s} \sum_{j} (\ga_i \gb_j \gd'_j + \gb_i \gb_j \gd'_j) \\
& \leq & \sum_{i \geq s} \sum_{j} (\ga_i \gb_j \gd'_i + \gb_i \gb_j \gd'_i) \\
& = &  \sum_{i \geq s} (\ga_i \gd'_i + \gb_i \gd'_i) \sum \gb_j  \\
& \leq &  \sum_{i \geq s} \ga_i \gc_i   \sum \gb_i  +
\sum_{i \geq s} \gb_i \gd'_i \sum \ga_i ,
\end{eqnarray*}
where we used:
$\sum \ga_i \gc_i \leq \sum \gb_i \gd'_i$ for the first inequality;
$\gd'_i \geq \gd'_j$ $\forall i \geq s$ (and $\forall j$) for the second;
and $\gc_i \geq \gd'_i$ $\forall i$ for the third.
This completes the proofs of Lemma \ref{ETS_for_aepi} and Theorem \ref{aepi}. \qed

\section{Counterexamples}\label{Ctrx}

Here we give the construction for Theorem \ref{Tctrex}.
As mentioned earlier, two further conjectures from \cite{Pem} turn
out to be disproved by the same examples, and we begin
by stating these.

\begin{conj}[\cite{Pem}, Conjecture 8]\label{conj:Pem8}
NA+ implies the NMP.
\end{conj}
(See Section \ref{section:ae} for NMP.)

Recall that for $\eta,\gz\in \Omega$, $\eta$ {\em covers} $\gz$ ($\eta\cov\gz$)
if there is an $i\in [n]$ for which $\eta_i=1$, $\gz_i = 0$ and
$\eta_j=\gz_j ~\forall j\neq i$.
Following \cite{Pem} we say that $\mu\in \m$
{\em stochastically covers} $\nu\in \m$ ($\mu\scov\nu$)
if we can couple random variables
$\eta $, $\gz$ having laws $\mu$ and $\nu$ so that with probability 1,
$\eta = \zeta$ or $\eta\cov\gz$; and that
$\mu$ has the {\em stochastic covering property} (SCP) if
$\mu(\cdot|\eta_i=0)\scov\mu(\cdot|\eta_i=1)$ for every $i$
(where, again, we regard these as measures on $\Omega_{[n]\sm \{i\}}$).
Observe that if $\mu$ is NA+ then
$\mu(\cdot|\eta_i=0) \succeq \mu(\cdot|\eta_i=1)$;
the following strengthening was suggested by Pemantle.
\begin{conj}[\cite{Pem}, Conjecture 9]\label{conj:Pem9}
NA+ implies the SCP.
\end{conj}
As already mentioned,
our examples invalidate both these conjectures.
\begin{thm}\label{thm:ctrex_Pem89}
Conjectures \ref{conj:Pem8} and \ref{conj:Pem9} are false.
\end{thm}
Conjecture \ref{conj:Pem8} was also disproved in \cite{BBL};
see also the note at the end of this section.

\medskip
We now describe the examples.
For a positive integer $k \geq 2$ and $\beta \in (0,1)$, let $\nu^{k,\beta}$ be the measure on $\Omega_{2k}$ with
$$\nu^{k,\beta}(\eta) \propto \left\{ \begin{array}{ll}
1 & \textrm{if } \ \ [|\eta|=k-1 \textrm{ and } \eta_1=1] \\
\beta^2 & \textrm{if } \ \ [|\eta|=k-1 \textrm{ and } \eta_1=0] \\
\beta & \textrm{if } \ \ |\eta|=k \\
\beta^2 & \textrm{if } \ \ [|\eta|=k+1 \textrm{ and } \eta_1=1] \\
1 & \textrm{if } \ \ [|\eta|=k+1 \textrm{ and } \eta_1=0] \\
0 & \textrm{otherwise}
\end{array} \right.$$
Note that (clearly) $\nu^{k,\beta}$ is almost exchangeable.
\begin{prop} \label{prop:ctrex}
The measure $\nu^{k,\beta}$ satisfies:
\\{\rm (a)} NA+ if and only if $\beta \geq \frac{1}{\sqrt{2}}$
\\{\rm (b)} ULC if and only if $\beta \geq 1-\frac{2}{k+1}$
\\{\rm (c)} unimodality (and also LC) if and only if
$\beta \geq 1-\sqrt{\frac{2}{k+1}}$
\\{\rm (d)} NMP if and only if $\beta \geq \sqrt{1-\frac{2}{k+1}}$
\\{\rm (e)} SCP if and only if $\beta \geq \sqrt{1-\frac{2}{k+1}}$
\end{prop}
For example, for $\gb = 0.71$:
$\nu^{6,\gb}$ is NA+ but not ULC, giving the first part of Theorem
\ref{Tctrex} (i.e. disproving Conjecture \ref{Pem4});
$\nu^{23,\gb}$ is NA+ but not unimodal (proving Theorem \ref{Tctrex}); and $\nu^{4,\gb}$ is NA+ but not NMP or SCP
(proving Theorem \ref{thm:ctrex_Pem89}).

\mn {\em Proof.}
We will mainly prove what we need for
Theorems \ref{Tctrex} and \ref{thm:ctrex_Pem89}, namely
``if" in (a) and ``only if" in (b)-(e).
The other direction in (b),(c) will come for free, but we omit
the (not very difficult) verifications of the remaining implications.

Fix $k$ and $\beta$, write $\nu$ for $\nu^{k,\beta}$, and set $r_i=\nu(|\eta|=i)$.
We have, for some $C$,
$$r_k= C \beta \C{2k}{k} \ \ \ \textrm{ and } \ \ \ r_{k-1}=r_{k+1}=C \Big( \frac{k-1}{2k}+\frac{k+1}{2k}\beta^2 \Big) \C{2k}{k-1}.$$
Unimodality and LC for $\nu$ are
equivalent to (each other and) $r_k \geq r_{k-1}$, which reduces to
$$\beta^2-2\beta+\frac{k-1}{k+1} \leq 0,$$
giving (c).
ULC for $\nu$ is equivalent to $k^2 r_k^2 \geq (k+1)^2 r_{k-1} r_{k+1}$, which reduces to
$$(k+1) \beta^2 - 2k \beta + (k-1) \leq 0,$$
giving (b).
The NMP requires that
$$\nu(\eta_1=1 \big| |\eta|=k)=\frac{1}{2}$$
be at least as large as
$$\nu(\eta_1=1 \big| |\eta|=k-1)=\frac{k-1}{(k-1)+(k+1)\beta^2},$$
from which the forward direction of (d) follows.
The SCP requires
$$\nu(|\eta|=k-1 \big| \eta_1=1) \leq \nu(|\eta|=k-1 \big| \eta_1=0),$$
which reduces to
$$\C{2k-1}{k-2} \leq \beta^2 \C{2k-1}{k-1}$$
and yields the forward direction of (e).

It remains to prove the backward direction of (a);
that is, we assume $\beta \geq 1/\sqrt{2}$ and should show $\nu$ is NA+.
Since $\nu$ is almost exchangeable, Theorem \ref{AE_equiv} says
we only need to show NC+, which, by symmetry, will follow if we show
$\eta_1 \da \eta_2$ and $\eta_2 \da \eta_3$
with respect to $W\circ \nu$, for any external field $W$.
(Our original proof of this has been shortened using some ideas from \cite{BBL}.)
Observe that, since a limit of NC measures is NC, it suffices to
consider the case when all entries of $W$ are finite and strictly positive.

Let $W'=(W_1,1 \dots 1)$, and let $\nu'$ be the projection of $W'\circ \nu$ on
$\Omega_{\{2 \dots 2k\}}$.
To prove $\eta_2 \da \eta_3$ for $W \circ \nu$, it suffices to show $\nu'$ is NC+,
which, since $\nu'$ is exchangeable, will follow {\em via}
Theorem \ref{E_equiv} if we show $\nu'$ has a ULC rank sequence.
The nonzero part of the normalized rank sequence
$(a_i:=\nu'(|\eta|=i)/\C{2k-1}{i})_{i=0}^{2k-1}$
is $(a_{k-2}\dots a_{k+1})\propto
(W_1,W_1\beta+\beta^2,W_1\beta^2+\beta,1)$, which a straightforward
computation shows to be LC when $\beta \geq 1/\sqrt{2}$.

That $\eta_1 \da \eta_2$ for $W \circ \nu$ will follow immediately from
\beq{SD1}
W \circ \nu (\cdot \big| \eta_1=0) \succeq W \circ \nu (\cdot \big| \eta_1=1).
\enq
Set
\begin{eqnarray*}
\pi_1 & = & W \circ \nu (\cdot \big| \eta_1=0,|\eta|=k+1), \\
\pi_2 & = & W \circ \nu (\cdot \big| \eta_1=0,|\eta| \in \{k-1,k\}), \\
\pi_3 & = & W \circ \nu (\cdot \big| \eta_1=1,|\eta| \in \{k,k+1\}), \textrm{ and} \\
\pi_4 & = & W \circ \nu (\cdot \big| \eta_1=1,|\eta|=k-1).
\end{eqnarray*}
It follows readily from Lemma \ref{rr_prod_nmp}
(since $\beta<1$ and the two measures appearing in (\ref{SD1}) are rank rescalings of
a common product measure,
namely the measure $\mu \in \m_{\{2 \dots 2k\}}$
with $\mu(\tau) \propto \prod W_i^{\tau_i}$) that each of $\pi_1$, $\pi_2$
stochastically dominates each of $\pi_3$, $\pi_4$.
Consequently, every convex combination of $\pi_1$ and $\pi_2$ stochastically
dominates every convex combination of $\pi_3$ and $\pi_4$, which in particular
gives (\ref{SD1}). \qed

Before closing this section, let us just mention that
a more natural class of counterexamples to
Conjecture \ref{conj:Pem8} is {\em probably} provided by the
following simple construction, which, as far as we know,
first appeared
in \cite{DG}.  Given $k$, let $G$ be the graph with
$V(G)= \{x,y,z_1\dots z_k\}$ and
$E(G) = \{xy,xz_1\dots xz_k,yz_1\dots yz_k\}$.
It is well known and easy to see
(consider the event $\eta_{xy}=1$)
that for $k\geq 5$, the
USF measure for $G$ fails the NMP, so is a counterexample to
Conjecture \ref{conj:Pem8} if the USF measure for $G$ is NA+.
The latter would follow from Conjecture \ref{fmconj}(b)
(USF measures are FM+) for $G$,
since Conjecture \ref{pef} for these graphs is contained in
the result from \cite{Wagner} mentioned following
Conjecture \ref{CRCM}.
(We can prove FM+ for $k\leq 5$,
and even this is not so easy).

\section{Proof of Theorem \ref{thm:rlc5}}\label{sec:rlc5}

The main point here is the following lemma, stating that NC+ implies the APP
(defined before Theorem \ref{TAPP}) for measures in $\m_4$.
For simplicity, we set $\ga_X=\mu(X)$ for $X \sub [n]$ (where we
now treat $\Omega_n$ as $2^{[n]}$),
often omit commas and set braces in subscripts
(e.g. $\ga_{134} = \mu(\{1,3,4\})$),
and write $\ga_0$ for $\mu(\emptyset)$.
Let $\gS_{r,s}^t = \sum \ga_X \ga_Y$, with the sum over unordered pairs
$\{X,Y\}$ of subsets of $[n]$ with $|X|=r$, $|Y|=s$, and $|X \cap Y|=t$.
\begin{lemma} \label{lemma:app4}
If $\mu \in \m_4$ is NC+, then
\beq{eq:app4}
3 \gS_{1,3}^0 \leq 4 \gS_{2,2}^0.
\enq
\end{lemma}
We first prove this and then give the easy derivation of Theorem \ref{thm:rlc5}.
(Notice that we could also
get Theorem \ref{thm:rlc5} from Lemma \ref{lemma:app4} via Theorem \ref{TAPP}, (\ref{eq:app4}) being
the only part of the CAPP that is not immediate from NC+;
but, as mentioned earlier, Theorem \ref{TAPP} is relatively difficult,
so we prefer a direct proof of Theorem \ref{thm:rlc5} here.)

For convenience, we now work with unnormalized
(nonnegative) measures on $\Omega$, and say that such a measure
$\mu$ with $\mu(\Omega)>0$
has a property (CNC, NC+, etc.) iff its normalization
$\mu'$
(given by $\mu'(\eta)=\mu(\eta)/\mu(\Omega)$) does.
Observe that $\A \da \B$ under
$\mu$
if and only if
\beq{eq:scale_inv_nc}
\mu(\A \bar{\B}) \mu(\bar{\A} \B) \geq \mu(\A \B) \mu(\bar{\A} \bar{\B})
\enq
(where $\bar{\A}=\Omega \sm \A$).

\mn{\em Proof of Lemma \ref{lemma:app4}.} Let $A = \gS_{2,2}^0$
($= \ga_{12}\ga_{34} + \ga_{13}\ga_{24} + \ga_{14}\ga_{23}).$
We may assume $\ga_0=\ga_{1234}=0$, since decreasing $\ga_0$ or $\ga_{1234}$ preserves NC+ and has no effect on (\ref{eq:app4}).
Furthermore, by renaming variables, applying a constant external field, and scaling (none of which affect (\ref{eq:app4})), we may assume
\beq{gaassume}
\ga_{123}=\ga_4=1 \ \ \textrm{ and } \ \ \ga_{1}\ga_{234},\ga_{2}\ga_{134},\ga_{3}\ga_{124} \leq 1.
\enq
(First, rename coordinates so $\ga_{i}\ga_{[4] \sm \{i\}}$ is largest for $i=4$ (and observe we may assume this largest value is strictly positive).
Second, impose a uniform external field $(W,W,W,W)$ to get $\ga_{123}=\ga_4$.
Third, divide all
values $\ga_X$ by $\ga_{123}$.)

Since $\mu$ is NC+ (so in particular CNC),
(\ref{gaassume}) gives
$\ga_{1} \leq \ga_{12}\ga_{13}$ (and similarly for $\ga_2$, $\ga_3$) and
$\ga_{124} \leq \ga_{14} \ga_{24}$ (and similarly for $\ga_{134}$, $\ga_{234}$).
It thus suffices to show
\beq{eq:lemmaets}
3(1 + xy + xz + yz) \leq 4A,
\enq
where
$$x=\ga_{12}\ga_{34}, \ y=\ga_{13}\ga_{24}, \textrm{ and } z=\ga_{14}\ga_{23} ~~ (\textrm{so } A=x+y+z).$$
For fixed $A$, the left hand side of (\ref{eq:lemmaets}) is maximized when $x=y=z$; thus (\ref{eq:lemmaets}) holds whenever $A \in [1,3]$ (as can be seen by examining the quadratic polynomial $A^2-4A+3$).
In view of (\ref{gaassume})
we can assume $A \leq 3$, so
we just need $A \geq 1$.

Assume, for a contradiction, that $A<1$.
Negative correlation of $\eta_2$ and $\eta_3$ for the measure
$(0,1,1,W) \circ \mu$ implies (use (\ref{eq:scale_inv_nc}))
$$P^1(W):=(\ga_{24}\ga_{34}-\ga_{234})W^2 + (\ga_2 \ga_{34} + \ga_3 \ga_{24} - \ga_{23})W + \ga_2 \ga_3 \geq 0 \ \ \ \ \textrm{ for } W>0.$$
Similarly, negative correlation of $\eta_1$ and $\eta_2$ for
$(\infty,1,1,W) \circ \mu$ implies
$$P_1(W):=\ga_{124}\ga_{134}W^2 +
(\ga_{12}\ga_{134}+\ga_{13}\ga_{124}-\ga_{14})W +
(\ga_{12}\ga_{13}-\ga_{1}) \geq 0 \ \ \ \ \textrm{ for } W>0.$$
Similarly (interchanging 1 with either 2 or 3) we have, again for $W>0$,
$$P^2(W):=(\ga_{14}\ga_{34}-\ga_{134})W^2 + (\ga_1 \ga_{34} + \ga_3 \ga_{14} - \ga_{13})W + \ga_1 \ga_3 \geq 0,$$
$$P_2(W):=\ga_{124}\ga_{234}W^2 + (\ga_{12}\ga_{234}+\ga_{23}\ga_{124}-\ga_{24})W + (\ga_{12}\ga_{23}-\ga_{2}) \geq 0,$$
$$P^3(W):=(\ga_{14}\ga_{24}-\ga_{124})W^2 + (\ga_1 \ga_{24} + \ga_2 \ga_{14} - \ga_{12})W + \ga_1 \ga_2 \geq 0, \textrm{ and}$$
$$P_3(W):=\ga_{134}\ga_{234}W^2 + (\ga_{13}\ga_{234}+\ga_{23}\ga_{134}-\ga_{34})W + (\ga_{13}\ga_{23}-\ga_{3}) \geq 0.$$

We pause to show
\beq{fullsupp}
\mbox{$\ga_X > 0$ for $X \neq \emptyset,[4]$.}
\enq
First we show $\ga_X > 0$ if $|X|=2$.
Suppose for example that $\ga_{12}=0$.
Since $P^1(W) \geq 0$ for all $W>0$ and $\ga_2=0$ (since $\ga_2\leq \ga_{12}\ga_{23}$)
the coefficient of $W$ in $P^1$ must be nonnegative, and thus
(using $\ga_3 \leq \ga_{13} \ga_{23}$)
$y\ga_{23} \geq \ga_{23}.$
If $\ga_{23} > 0$, this gives
$A \geq y \geq 1$;
thus (since we are assuming $A<1$) $\ga_{23}=0$, and similar reasoning shows $\ga_{1}=\ga_{3}=\ga_{13}=0$.
Hence, $\ga_X=0$ unless $X=\{1,2,3\}$ or $4 \in X$; but then nonnegativity of $P_1$, $P_2$, and $P_3$ gives $\ga_{14}=\ga_{24}=\ga_{34}=0$.
Thus $\ga_X>0$ if and only if $X=\{1,2,3\}$ or $X=\{4\}$; but for any such
measure $\eta_1$ and $\eta_2$ are strictly positively correlated.
This contradiction shows $\ga_{12}>0$, and
similar arguments (or symmetry) give $\ga_X > 0$
whenever $|X|=2$.

If $\ga_2=0$, then nonnegativity of the linear term in $P^1$ gives, as
in the preceding
paragraph
(and using $\ga_{23}>0$), $A \geq 1$; thus $\ga_2>0$.
Similar arguments (or, again, symmetry)
show $\ga_1$, $\ga_3$, $\ga_{124}$, $\ga_{134}$, and $\ga_{234}$ are positive, and we have (\ref{fullsupp}).

Set
$$a=\frac{\ga_1}{\ga_{12}\ga_{13}}, \ b=\frac{\ga_2}{\ga_{12}\ga_{23}}, \ c=\frac{\ga_{3}}{\ga_{13}\ga_{23}}, \ d=\frac{\ga_{124}}{\ga_{14}\ga_{24}}, \ e=\frac{\ga_{134}}{\ga_{14}\ga_{34}}, \ \textrm{ and } \ f=\frac{\ga_{234}}{\ga_{24}\ga_{34}}.$$
Note $a,b,c,d,e,f \in (0,1]$.
If the coefficient of W in $P^1$ is nonnegative, then, as above, $A \geq 1$; thus this coefficient is negative, whence the discriminant of $P^1$ is nonpositive.
This yields
$$1-\frac{\ga_2 \ga_{34}}{\ga_{23}}-\frac{\ga_3 \ga_{24}}{\ga_{23}} \leq 2 \sqrt{\frac{\ga_2 \ga_3}{\ga_{23}^2}(\ga_{24}\ga_{34}-\ga_{234})},$$
which in the notation introduced above becomes
$$1-bx-cy \leq 2 \sqrt{bxcy(1-f)} \leq (bx+cy) \sqrt{1-f}.$$
Thus
$$x+y \geq \big[ (1+\sqrt{1-f}) \max{\{b,c\}} \big]^{-1},$$
and a similar argument using $P_1$ gives
$$x+y \geq \big[ (1+\sqrt{1-a}) \max{\{d,e\}} \big]^{-1},$$
so that
$$x+y \geq \max{\big\{ \big[ (1+\sqrt{1-f}) \max{\{b,c\}} \big]^{-1},\big[ (1+\sqrt{1-a}) \max{\{d,e\}} \big]^{-1} \big\}}.$$
Similar arguments using $P^2$, $P_2$, $P^3$, and $P_3$ yield
$$x+z \geq \max{\big\{ \big[ (1+\sqrt{1-e}) \max{\{a,c\}} \big]^{-1},\big[ (1+\sqrt{1-b}) \max{\{d,f\}} \big]^{-1} \big\}}$$
and
$$y+z \geq \max{\big\{ \big[ (1+\sqrt{1-d}) \max{\{a,b\}} \big]^{-1},\big[ (1+\sqrt{1-c}) \max{\{e,f\}} \big]^{-1} \big\}}.$$
In particular, we have $x+y,x+z,y+z \geq 1/2$.

The proof is now an easy consequence of
\beq{calc_claim}
\inf \big\{ \big[(1+\sqrt{1-v})u \big]^{-1} +
\big[(1+\sqrt{1-u})v \big]^{-1} : 0 < u,v \leq 1 \big\}=\frac{27}{16},
\enq
verification of which is a straightforward calculus exercise which we omit.
Assuming (\ref{calc_claim}) and, without loss of generality,
$a \geq b$ and $d \geq e$, we have
$$(x+y)+(y+z) \geq \big[ (1+\sqrt{1-a})d \big]^{-1} +
\big[ (1+\sqrt{1-d})a \big]^{-1} \geq \frac{27}{16},$$
which, combined with $x+z \geq 1/2$, gives the final contradiction $2A > 2$. \qed

\mn{\em Proof of Theorem \ref{thm:rlc5}.} First observe that LC[$m$] is equivalent to having
\beq{eq:rulc}
\mbox{NC+ implies ULC for measures in $\m_n$}
\enq
for all $n \leq m$.
Suppose $\mu \in \m_{n}$ has rank sequence $(r_i)_{i=0}^n$.
Notice that in general for (\ref{eq:rulc}) it is enough to show that NC+ implies
\beq{eq:muulck}
r_k^2 \C{n}{k}^{-2} \geq r_{k-1} r_{k+1} \C{n}{k-1}^{-1} \C{n}{k+1}^{-1}
\enq
for $1 \leq k \leq \lfloor n/2 \rfloor$, since the measure $\mu^{*} \in \m_{n}$ with $\mu^{*}(X) = \mu([n] \sm X)$ has rank sequence $(r_{n-i})_{i=0}^n$ and is NC+ if and only if $\mu$ is.
In fact, Choe and Wagner \cite{choewagner03} show that (\ref{eq:muulck}) holds
for $k=1$ and {\em any} $n$ (assuming NC+). This gives (\ref{eq:rulc}) for
$n \leq 3$ (and hence LC[3]) and for the cases of interest here---that is,
$n=4,5$---reduces the problem to proving (\ref{eq:muulck}) when $k=2$.

Assume $n \in \{4,5\}$.
Using inequalities of the form
$$\ga_{i} \ga_{ijl} \leq \ga_{ij} \ga_{il}$$
(which follow from NC+) we obtain
\beq{eq:sigma1}
\gS_{1,3}^1 \leq \gS_{2,2}^1.
\enq
It follows from Lemma \ref{lemma:app4} (for $n=4$ this is the conclusion of the lemma, and for $n=5$ we apply the lemma to each of the five conditional measures $\mu(\cdot | \eta_i=0)$) that
\beq{eq:sigma0}
3\gS_{1,3}^0 \leq 4 \gS_{2,2}^0.
\enq
Note also that Cauchy-Schwarz implies that the average size of a term in
$\gS_{2,2}^2$ is at least
the average size of a term in
either of $\gS_{2,2}^0$, $\gS_{2,2}^1$; that is,

\beq{eq:sigma2}
\gS_{2,2}^2 \geq \frac{4}{(n-2)(n-3)} \gS_{2,2}^0
\ \ \ \textrm{ and } \ \ \
\gS_{2,2}^2 \geq \frac{1}{n-2} \gS_{2,2}^1.
\enq
Thus, finally, we have (\ref{eq:muulck}) for $k=2$:
$$9 r_1 r_3 = 9 \gS_{1,3}^0 + 9 \gS_{1,3}^1
\leq 12 \gS_{2,2}^0 + 9 \gS_{2,2}^1
\leq 8 \gS_{2,2}^0 + 8 \gS_{2,2}^1 + 4 \gS_{2,2}^2
= 4 r_2^2 \ \ \ \textrm{ if } n=4$$
and
$$2 r_1 r_3 = 2 \gS_{1,3}^0 + 2 \gS_{1,3}^1
\leq \frac{8}{3} \gS_{2,2}^0 + 2 \gS_{2,2}^1
\leq 2 \gS_{2,2}^0 + 2 \gS_{2,2}^1 + \gS_{2,2}^2
= r_2^2 \ \ \ \textrm{ if } n=5,$$
where in each case we used (\ref{eq:sigma1}) and (\ref{eq:sigma0}) for the first inequality and (\ref{eq:sigma2}) for the second. \qed

\section{Urns}\label{Urns}

Finally, in this short section,
we just give the easy examples justifying
Proposition \ref{Uctrex} and the remark following it.
(Recall that these say that log-concavity and the Rayleigh
property fail for competing urn measures
(with identical balls).
As mentioned in Section \ref{Intro}, positive results for
competing urns will appear separately.)
In both examples we use $p(j)$ for the probability
that any given ball lands in urn $j$.

\begin{ex}\label{Uex1}
Suppose we have three balls and urns $0\dots n$,
with $p(0)=\eps$ and $p(1)=\cdots = p(n)=(1-\eps)/n$, where
$\eps$ is small and
$n\eps^{3/2}$ is large.
Then for the associated rank sequence, say $a = (a_1,a_2,a_3)$,
we have $a_1\approx\eps^3$, $a_3\approx (1+2\eps)(1-\eps)^2$ and
$$a_2 = 3\eps^2(1-\eps) + 3\eps(1-\eps)^2/n + 3(1-\eps)^3(n-1)/n^2
\approx 3\eps^2(1-\eps);$$
so LC fails for a.
\end{ex}
(We don't know
what happens if we replace ``LC" by ``unimodal.")

\begin{ex}\label{Uex2}
Suppose we have two balls, urns $0,1,2$, and
$p(1)=p(2)=\eps$, with
$\eps$ small,
and impose the external field $(\eps,1,1)$.
Then for the corresponding urn measure $\mu$ on $\{0,1\}^{\{0,1,2\}}$
(and $\eta$ the random configuration) we have
$\mu(\eta_1=\eta_2 =1) \propto 2\eps^2$,
$\mu(\eta_1=\eta_2 =0) \propto (1-2\eps)^2\eps$, and
$\mu(\eta_1=1, \eta_2 =0) ,\mu(\eta_1=0, \eta_2 =1)
\propto\eps^2 +2\eps^2(1- 2\eps)$,
so that $\eta_1$ and $\eta_2$ are strictly positively correlated.
\end{ex}

\bn
{\bf Acknowledgments.}  Thanks to David Wagner for helpful comments and
some pointers to the literature,
and thanks to the referees for carefully reading the paper.

\end{document}